\documentclass[12pt]{article}
 \usepackage[english,activeacute]{babel}
 \usepackage[pdftex]{graphicx}
 \usepackage[latin1]{inputenc}
 \pdfoutput=1
 \usepackage{amsmath,color,graphics}
 \usepackage{hyperref,amssymb}
 \usepackage{ae}
 
 \title{AN ALGORITHM FOR SINGULAR VALUE DECOMPOSITION OF MATRICES IN BLOCKS\\
      \normalsize{Technical Report}}
 \author{\'Avaro Francisco Huertas-Rosero}
\begin{document}
 \maketitle
 \begin{abstract}
  Two methods to decompose block matrices analogous to Singular Matrix Decomposition are proposed, one yielding the so called economy decomposition, and other yielding the full decomposition.  This method is devised to avoid handling matrices bigger than the biggest blocks, so it is particularly appropriate when a limitation on the size of matrices exists.   The method is tested on a document-term matrix (17780$\times$3204) divided in 4 blocks, the upper-left corner being 215$\times$215.
 \end{abstract}

 \section{Introduction}
  Singular Value Decomposition has proved to be useful in a wide range of applications, where a linear relation is a suitable model for a big number of variables.   Its main strength is in its ability to abstract most of the meaningfull relation in a much smaller subspace~\cite{StructuresLSA},\cite{SVDProjections},\cite{LSA},\cite{Latent}.

  Even though the calculations are very simple in essence, the method is at its best when dealing with big dimension matrices of data, and the computational resources to perform the calculations are often insufficient.

  In this document I propose an algorithm that allows to deal with the matrix by pieces, so it does not need to define big matrices or operate with them, but only smaller blocks.  
 \section{The usual algorithm}the usual algorithm to perform the decomposition is made in two steps: First, a transformation is found that takes the matrix to bidiagonal form, and then, the bidiagonal matrix is decomposed with a different procedure.
  \subsection{Householder transformations and bidiagonal matrices}
   The first step is carried by means as a certain class of symmetrical orthogonal (or unitary) matrices called Householder transformations~\cite{Householder}.   A Householder transformation is defined by a unitary vector this way:
   \begin{equation}
    H_R = 1_{n\times n} - 2(\hat{R})^t\hat{R}
   \end{equation}
   where $\hat{R}$ is an unitary vector of dimension (number of components) n, and $1_{n\times n}$ is the identity with n rows and columns.   It is easy to see that the matrix that corresponds to this transformation is symmetric, which means that it is not changed by transposition (changing rows by columns, and vice versa).  That impplies that it is its own inverse, e.i., that its square is the identity matrix.

   A householder operation can be found that, when multiplied by the left, turns the all but one of the entries of the first column of a matrix into zero, but preserving the sum of the squares of the entries of that column.
   \begin{equation}
    \left(1_{n\times n} - 2(\hat{R})^t\hat{R}\right)
    \begin{pmatrix}
     a\\
     d\\
     g\\
     k
    \end{pmatrix} =
    \begin{pmatrix}
     \sqrt{a^2 + d^2 + g^2 + k^2}\\ 0 \\ 0 \\ 0
    \end{pmatrix}
   \end{equation}

   The unitary vector that defines the Householder transformation can be computed as having one part proportional to the part that is to turn into zero.  Being an unitary vector, the proportionality factor can be best represented by some unknown factor $X$ divided by the norm of that part of the vector, the square root of the sum of the squares.  In our example:
   \begin{equation}
    \hat{R} = \left(\sqrt{1 - X^2}, \frac{Xd}{\sqrt{d^2 + g^2 + k^2}},\frac{Xg}{\sqrt{d^2 + g^2 + k^2}},\frac{Xk}{\sqrt{d^2 + g^2 + k^2}}k\right)
   \end{equation}

   Imposing the condition that it makes the required entries of the vector 0, the unknown factor turns out to be:
   \begin{equation}
    X = \frac{a + \sqrt{a^2 + d^2 + g^2 + k^2}}{2\sqrt{a^2 + d^2 + g^2 + k^2}} = \frac{a + Norm}{2Norm}
   \end{equation}

   We can express any matrix in block form, separating the first row and column, in the same way as the Householder matrix.  This latter has a very simple form:

   \begin{multline}
    M = \begin{pmatrix}
     M_0 & M_{row}\hat{V}_{row}\\
     M_{col}(\hat{V}_{col})^t & M_{block}
    \end{pmatrix} \\
    H_M = \frac{1}{N}
    \begin{pmatrix}
     M_0 & M_{col}\hat{V}_{col}\\
     M_{col}(\hat{V}_{col})^t & -N\cdot 1_{(n-1)\times (n-1)}+(N-M_0)\hat{V}^t\hat{V}
    \end{pmatrix}
   \end{multline}
   where $M_0$ is the first diagonal entry of the matrix, and the remaining of the first row and first column are expressed by a norm and an unit vector, $M_{row}\hat{V}_{row}$ and $M_{col}(\hat{V}_{col})^t$.   N would be the norm of the first column, that is $N = \sqrt{M_0^2+M_{col}^2}$
   
   Multiplying the matrix $M$ by the Householder matrix by left, we get the following result:
   \begin{equation}
    H_M M = \frac{1}{N}
    \begin{pmatrix}
     N^2 &  - M_{col}M_{row}\hat{V}_{row} + M_{col}\hat{V}_{col}M_{block}\\
     0 &  M_{col}M_{row}(\hat{M}_{col})^tM_{row} - NM_{block} + (N - M_0)(\hat{M}_{col})^t\hat{M}_{col}M_{block}
    \end{pmatrix}
   \end{equation}

   A Householder transformation applied on the right, can have the same effect of annihilating all but one of the entries of the first row.  But, if we try to get a diagonal block matrix using two consecutive Householder transformation, one on the left and one on the right, the second is going to turn the zero entries produced by the first into other number, thus failing to produce a block diagonal matrix.

  \begin{center}
   \includegraphics{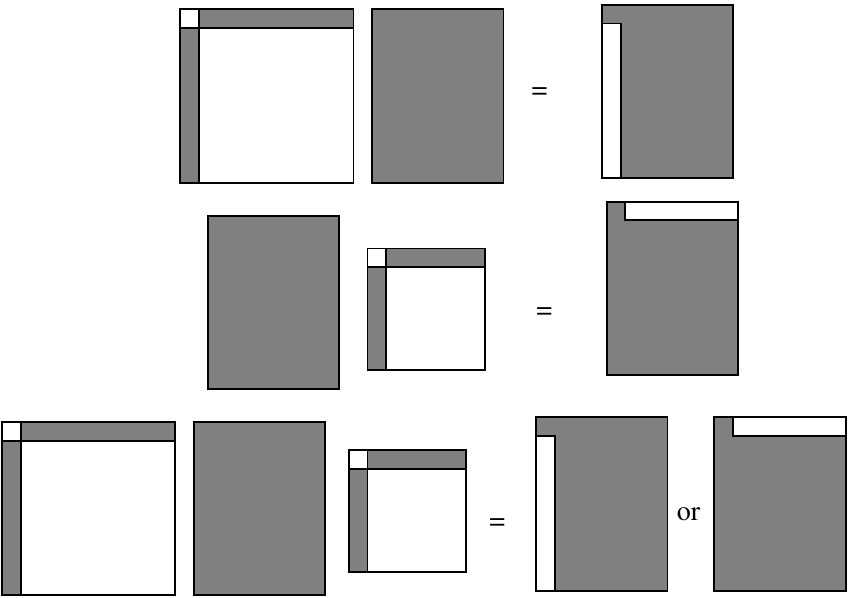}
  \end{center}

   We can, however, leave one nonzero element in the first column, and two nonzero elements in the first row, using on the left a Householder transformation that mixes the contents of the rows only from the second element on:
  \begin{center}
   \includegraphics{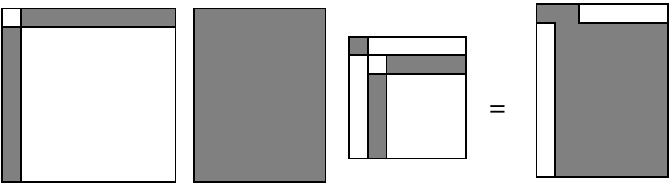}
  \end{center}

  Next we can focus on the right-bottom block of the matrix where the zero elements have not yet been produced.  The right and left Householder transformation for that block can be computed, so it will be left with only two nonzero elements in the first row and column.    Then, the same is done for the remaining block, and so on.
  \begin{center}
   \includegraphics{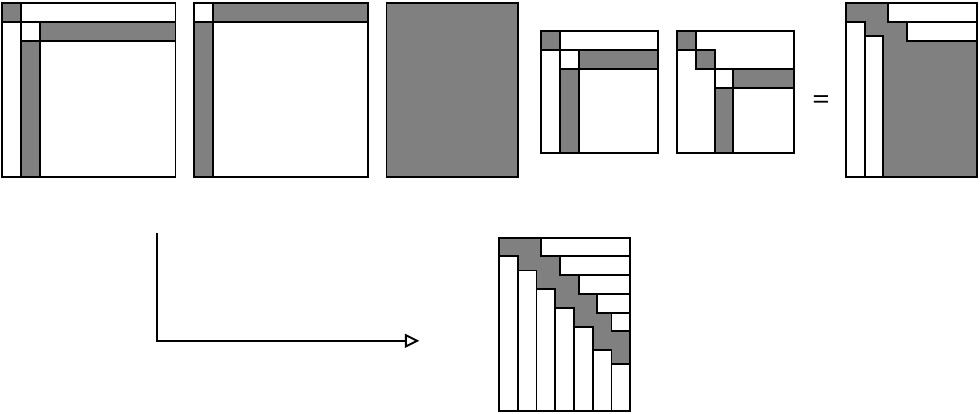}
  \end{center}
  The number of steps necessary to bring the matrix to a bidiagonal form has been $2n-1$, where $n$ is the lower dimension of the matrix.  It is a very fast step.

  \subsection{Iterative diagonalisation}
   To take the bidiagonal matrix to a diagonal form is not so easy, and cannot be done in a fixed number of steps, but has to be an iterative process.   The most efficient method for this is the \textbf{QR factorisation}~\cite{QRdecomposition}.

   Starting from our upper bidiagonal matrix, the steps to be performed are:
   \begin{enumerate}
    \item Construct the orthogonal transformation to bring the matrix to \textit{lower} triangular form, with Householder transformations like the ones we used in the last part, but only applied by right.
    \begin{center}
     \includegraphics{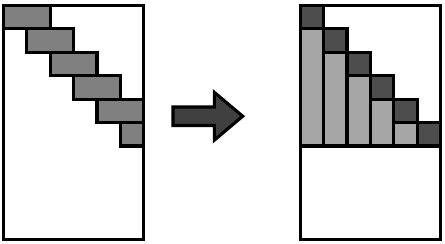}
    \end{center}
    \item Apply the transformation by left as well.  The result will be a matrix that is not either upper or lower triangular, but has the values more concentrated on the diagonal.
    \begin{center}
     \includegraphics{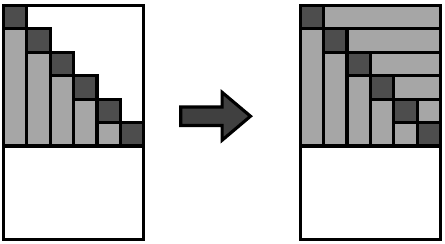}
    \end{center}
    \item Repeat the procedure.
    \begin{center}
     \includegraphics{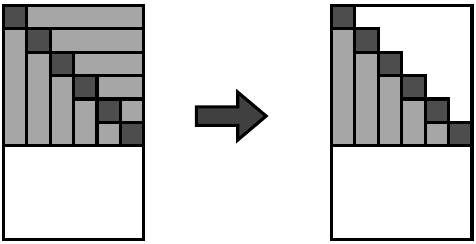}
    \end{center}
   \end{enumerate}

   All the steps are done with orthogonal transformation, and that ensures that the singular values are preserved.  In each step, besides, the square of every diagonal element are increased with the squares of the other elements in the row.   That, together with the preservation of singular values, ensures convergence.

  \subsection{getting the economy decomposition}
   The last procedure yields the two square unitary matrices, and the singular values diagonal must be nonsquare.   That means that a lot of memory is needed for the unitaries, wich are huge matrices with double precision.   A lighter alternative is to compute the \textit{economy} decomposition~\cite{EconomySVD}, wich give us just slices of unitary matrices.  The left one has as much rows as our matrix, but only as much columns as the rank (number of nonzero singular values) of the matrix.   And the right one has as much columns as our matrix, but only as much rows as the rank of the matrix.

   \begin{multline}
    M = U D V^t \\
    U^tU = 1_{rank} \hspace{2cm} V^tV = 1_{rank} \\
    UU^t \neq 1_{rows} \hspace{1cm} VV^t \neq 1_{columns}
   \end{multline}
   \begin{center}
    \includegraphics{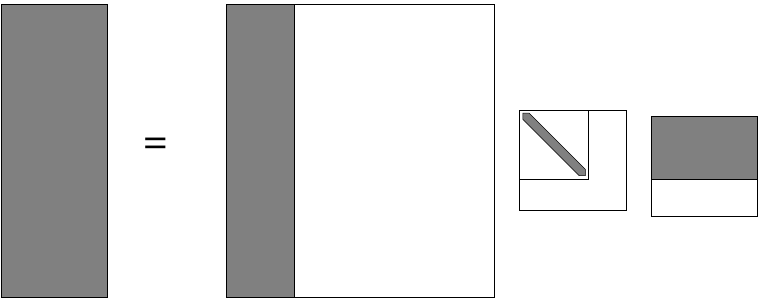}
   \end{center}

   To compute the economy decomposition, a lesser computational effort is needed.   Suppose that we have in our matrix more rows than columns, or that we took it to that form by transposing it.   Then, the symmetric matrix $A^tA$ will be smaller in dimensions than A.  We can use one of the usual algorithms to diagonalise it, and get
   \begin{equation}
    A^tA = V D V^t
   \end{equation}
   It is very easy to compute the inverse of square root of this matrix, all that is necessary is to take inverses of square roots of the diagonal elements of D.   Then, we can express our matrix like this:
   \begin{equation}
    A = A (A^tA)^{-1/2}(A^tA)^{+1/2} = \left(A V D^{-1/2}\right) D^{+1/2} V^t = U D^{+1/2} V^t
   \end{equation}
   It can be easily shown that $U=A V D^{-1/2}$ is a slice of a unitary matrix:
   \begin{equation}
    U^tU =  \left( D^{-1/2} V^t A^t\right)\left(A V D^{-1/2}\right) = D^{-1/2} V^t \left(A^tA\right) V D^{-1/2} = D^{-1/2} D D^{-1/2} = 1
   \end{equation}
   If the square matrix $A^tA$ hapens to be singular, then $V$ is also a slice of a square unitary, and $D$ is smaller, but allways invertible.

 \section{Frobenius norm and a better starting point}
  Latent Semantic Analysis, and other techniques, are based on the fact that some big matrices can be accurately represented only by the bigger terms of their spectral decomposition, that is, only the bigger singular values.    The usual convention to represent the diagonal matrix of singular values is in an ordered form, from the biggest, in the first element, to the smallest.

  The matrices are usually rather sparse, and with only some column and row swaps, we can take the bigger elements up and to the left, so the matrix is going to be closer to the desired form.

  The singular value problem, as the eigenvalue problem, can be seen as a maximisation of a certain value.  The solution of the following problem gives the left and right eigenvectors corresponding to the highest singular value:
 \begin{align}
  \text{LEFT: }maximise\ \langle \Psi \vert AA^t \vert \Psi \rangle &\ constrained\  by\ \langle \Psi \vert \vert \Psi \rangle = 1\\
  \text{RIGHT: }maximise\ \langle \Phi \vert A^tA \vert \Phi \rangle &\ constrained\  by\ \langle \Phi \vert \vert \Phi \rangle = 1
 \end{align}
  All along this work, the maximisation of the values in the upper left corner of the matrix is going to be used to arrive near to the diagonal form.  The first thing that can be done, is just arranging the rows and columns to take the higher values to the upper left corner of the matrix.

  There is a well known result, that tells us that the Frobenius norm of a matrix~\cite{Frobenius}, that is, the trace of its square $\Vert A \Vert = \sqrt{Trace(A^tA)}$, is the sum of the squares of the singular values, and is as well the sum of the squares of all the elements of the matrix.
  \begin{equation}
   Trace(A^tA) = \sum_i\left(\sum_j (A^t)_{ij}A_{ji}\right)= \sum_{ij}(A_{ij})^2
  \end{equation}
  If we consider the other square matrix $AA^t$ the result is the same, because in that case we just swap indexes $i$ and $j$.    With the SVD decomposition of the matrix, we only need to remember that a unitary matrix does not affect the trace:
  \begin{align}
   Trace(A^tA) &= Trace\left(\left(V (D_{mn})^t U^t\right)(U D_{mn} V^t) = Trace((D_{mn})^tD_{mn})\right) = \sum_i (\lambda_{ii})^2\\
  Trace(AA^t) &= Trace\left(\left(U (D_{mn})^t V^t\right)(V D_{mn} U^t) = Trace((D_{nm})^tD_{nm})\right) = \sum_i (\lambda_{ii})^2
  \end{align}
  where $D_{mn}$ is a square $m\times n$ diagonal matrix with the singular values. 

  A vector can be computed with the norms of each row, and the Cartesian norm of this vector will be the frobenius norm of the matrix.   The same can be done with the columns.   This two vectors can be used to sort the rows and columns of the matrix to get the higher values in the upper left corner.

 \begin{center}\includegraphics{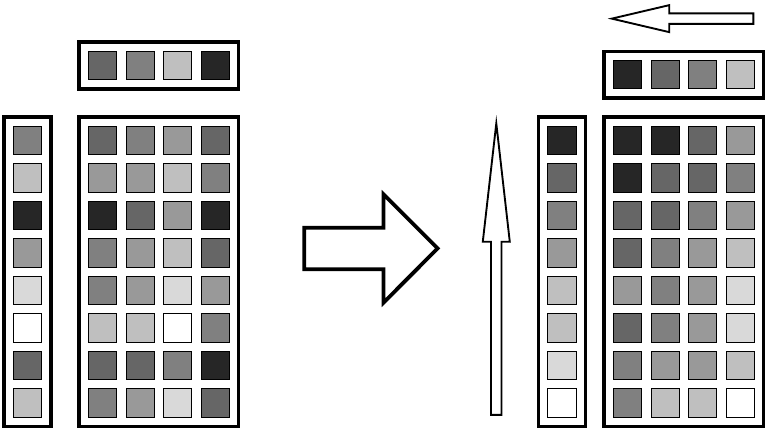}\end{center}

  Afther that sorting, a definition of the blocks can be done, with some criterion based on those row and column norms.    The blocks can be defined, for example, as to put a certain percentage of the whole frobenius norm in the first column-block, and a certain percentage on the first row-block. 

  If only the higher singular values are needed, it is not necessary to decompose the whole matrix, but instead two steps can be taken:
  \begin{enumerate}
   \item Separate the subspace of the highest singular values from that of the lower singular values
   \item Decompose only the block corresponding to the highest singular values
  \end{enumerate}
 \section{Partial SVD}
  In the proposed algorithm, the matrix is prepared so as to have more rows than columns (transposed if necessary) and the columns are cut in such a way that a fraction of the total square Frobenius norm is enclosed in the first column block.   The rows are separated in such a way to have square upper left block.

  Then, each block can be considered separately, and that can require considerably less computer resources.

 \section{A generalisation of Householder matrices for blocks}
  The first thing that can be done, is to generalise the concept of a Householder transformations to any partition of the rows and columns of the matrix in four blocks.  The general form of such a transformation is the following:
  \begin{equation}
   H = \begin{pmatrix}\label{DefinitionHH}
    1_{n\times n} - U(1_{m\times m}-\alpha)U^t & U \beta V^t \\
    V \beta U^t & 1_{(N-n)\times(N-n)}  V(1_{m\times m}+\alpha)V^t
   \end{pmatrix} \hspace{1cm} 
  \end{equation}
  where $U$ is a slice of a unitary matrix with $n$ rows and $m$ columns, and $V$ is a slice of unitary matrix with $N-n$ rows and $m$ columns.   This transformation would be a $N\times N$ unitary matrix, to be multiplied by left to a matrix with $N$ rows divided in two blocks with $n$ and $N-n$ rows, or by right to a matrix with $N$ columns split in blocks with $n$ and $N-n$ columns.   $\alpha$ and $\beta$ are diagonal matrices of rank $m$ with the property $\alpha^2+\beta^2 = 1_{m\times m}$

  This matrix y symmetric and it is its own inverse, two of the properties of a Householder transformation.   A householder transformation shifts the sign of only one vector, but this transformation can be shown to change the sign of any vector that lies within a subspace.  This subspace is defined by a set of mutually orthogonal vectors, which can be arranged in a column or row block.   The matrix can be written also like this:

  \begin{equation}
    H = 1_{N\times N} - 2\left(\begin{pmatrix}
	  U\sqrt{\frac{1}{2}(1-\alpha)} \\
     -V\sqrt{\frac{1}{2}(1+\alpha)}
    \end{pmatrix}
    \begin{pmatrix}
	  \sqrt{\frac{1}{2}(1-\alpha)}U^t &
	  -\sqrt{\frac{1}{2}(1+\alpha)}V^t
    \end{pmatrix}\right)
  \end{equation}
  A transformation like this can be used, for example, to annihilate a block of a matrix, just as in the usual SVD method.  Here are the steps to annihilate a nondiagonal block by multiplication by left:
  \begin{enumerate}
   \item The two relevant blocks, the ones that are going to be transformed to annihilate one of them, are decomposed.   Full SVD is not necessary, a simple decomposition Unitary-Symmetric will do.
   \begin{equation}
    A_{ij} = U_{ij}S_{ij} \hspace{2cm} S_{ij} = \left((A_{ij})^tA_{ij})\right)^{1/2} \hspace{2cm} U_{ij} = A_{ij}S_{ij}^{-1}
   \end{equation}
   \item The unitary factors of the two blocks are taken as the $U$ and $V$ matrices of the Householder matrix
   \item If we take the first column of a 2x2 blocks matrix as the relevant blocks, the action of the Householder matrix will give in the nondiagonal block:
   \begin{equation}
    \begin{pmatrix}
    1 - U(1-\alpha)U^t & U \beta V^t \\
    V \beta U^t & 1-  V(1+\alpha)V^t
    \end{pmatrix}
    \begin{pmatrix}
    US_{11} & A_{12} \\
    VS_{21} & A_{22}
    \end{pmatrix}=
    \begin{pmatrix}
     \#\#\# & \#\#\# \\
    V(\beta S_{11} -\alpha S_{21}) & \#\#\#
    \end{pmatrix}
   \end{equation}
   For this block to be zero, the parameters of the matrix must be:
   \begin{equation}
    \alpha = (1_{m\times m} + S_{21}(S_{11})^2S_{21})^{-1/2}\hspace{2cm}
    \beta = (1_{m\times m} - \alpha^2)^{1/2}
   \end{equation}
  \end{enumerate}
  There is another way of doing it as well, wich will probably take more time, but is based on a well known technique: the \textbf{GSVD: Generalised Single Value Decomposition}.   This is the simultaneous decomposition of two matrices:
   \begin{equation}
    A_{11} = UCX^t\hspace{2cm}
    A_{21} = VSX^t
   \end{equation}
  where, in the economy representation, $U$ and $V$ are slices of unitary matrices with the same dimensions than those of the other method.  X is a square matrix.   The matrices $S$ and $C$ are square and diagonal, and fulfil $C^2+S^2=1$.  They can be used indeed as $\alpha$ and $\beta$ respectively ($C=\alpha,S=\beta$) in the Householder matrix.

  With this kind of Householder transformations, we can perform a complete (not economy) blockwise decomposition of a matrix, iterating the annihilation of the two nondiagonal blocks as shown in the figure:

  \section{Blocks and the trace trick}
  To be able to perform the decomposition directly in the economy representation, a version of the \textbf{eigenvalue} (spectral) decomposition is needed.   The formula shown for annihilating blocks only works multiplied by one side, but it does not work to annihilate nondiagonal blocks acting on both sides, as an equivalence transformation.

  For a 2x2 number symmetric matrix, the problem of diagonalising it amounts to finding a certain number $x$ that fullfills:
  \begin{equation}
   \begin{aligned}
   U M U^t=&
   \begin{pmatrix} \sqrt{1-x^2} & x \\ x & -\sqrt{1-x^2} \end{pmatrix}
   \begin{pmatrix} A & B \\ B & C \end{pmatrix}
   \begin{pmatrix} \sqrt{1-x^2} & x \\ x & -\sqrt{1-x^2} \end{pmatrix}
  = D \\
    D &= \begin{pmatrix}
    \frac{1}{2}(A+B)+\sqrt{\frac{1}{2}(A-B)^2 + C^2} & 0 \\
    0 & \frac{1}{2}(A+B)-\sqrt{{\frac{1}{2}(A-B)^2 + C^2}}
   \end{pmatrix}
   \end{aligned}
  \end{equation}
  The condition is better derived from the null elements of the matrix, and is:
  \begin{equation}
   (1-x^2)B + x\sqrt{1-x^2}A =  x\sqrt{1-x^2}C + x^2 B
  \end{equation}
  Defining $\alpha = (1-2x^2)$ the condition becomes very simple, because $\sqrt{1-\alpha^2} = 2x\sqrt{1-x^2}$
  \begin{equation}
   \alpha B =  \frac{1}{2}\sqrt{1-\alpha^2}(A-C)
  \end{equation}
 The solution is easily found to be:
  \begin{equation}
   x = \sqrt{\frac{\sqrt{(A-C)^2 + 4B^2}-(A-C)}{2\sqrt{(A-C)^2 + 4B^2}}}
  \end{equation}
 On the other hand, if the entries of the matrix are suitably sized blocks, the condition is a lot more complicated.  We can represent the unitary matrix as being constructed with blocks $x$, $y$ and $z$ having the form shown above for Householder matrices.
  \begin{equation}
   \begin{aligned}
   U M U^t=&
   \begin{pmatrix} x & y \\ y^t  & z \end{pmatrix}
   \begin{pmatrix} A & B \\ B & C \end{pmatrix}
   \begin{pmatrix} x & y \\ y^t  & z \end{pmatrix}\\
    &=
   \begin{pmatrix} xAx + yB^tx + xBy^t + yCy^t & xAy + yB^ty + xBz + yCz\\
    y^tAx + zB^tx + y^tBy^t + zCy^t &  y^tAy + zB^ty + zCz
    \end{pmatrix}
   \end{aligned}
  \end{equation}
  The noncommutativity of the matrices does not allow for an easy solution as that with numbers.  Furthermore, to solve the condition for the nondiagonal elements should not be possible, except in the 2x2 or 2x3 blocks case, because that could be translated to solve analytically a general equation of order higher than five.  That, according to Abel's theorem, is not possible.

  But there is something we can do, and it is working with traces.  We can either maximise the trace of the first diagonal block, or minimise the trace of the square of the nondiagonal block.

  On the other hand, using slices of the unitary matrices, and their complement (the slice that is lacking for the total unitary) we can build an unitary matrix $S$ that allow us to isolate just a subspace to work on it, thus reducing substantially the dimensionality of the problem.
   \begin{equation}\label{UnitarySubspaces}
    S = \begin{pmatrix} U_1 & \bar{U}_1 & 0 & 0 \\0 &  0 & U_2 & \bar{U}_2 \end{pmatrix}
   \end{equation}
   Note that if block $(A^tA)_{12}$ is full rank, then $U_1$ would be square and there would not be a $\bar{U}_1$.

  Applying this unitary transformation to the matrix we get:
   \begin{equation}
    S^t (A^tA) S =
   \begin{pmatrix}
   (U_1)^t (A^tA)_{11} U_1 & (U_1)^t (A^tA)_{12} \bar{U}_1 & (U_1)^t (A^tA)_{11} U_2 & (U_1)^t (A^tA)_{12} \bar{U}_2 \\
   (\bar{U}_1)^t (A^tA)_{21} U_1 & (\bar{U}_1)^t (A^tA)_{22} \bar{U}_1 & (\bar{U}_1)^t (A^tA)_{21} U_2 & (\bar{U}_1)^t (A^tA)_{22} \bar{U}_2 \\
   (U_2)^t (A^tA)_{11} U_1 & (U_2)^t (A^tA)_{12} \bar{U}_1 & (U_2)^t (A^tA)_{11} U_2 & (U_2)^t (A^tA)_{12} \bar{U}_2 \\ 
   (\bar{U}_2)^t (A^tA)_{21} U_1 & (\bar{U}_2)^t (A^tA)_{22} \bar{U}_1 & (\bar{U}_2)^t (A^tA)_{21} U_2 & (\bar{U}_2)^t (A^tA)_{22} \bar{U}_2 
    \end{pmatrix}
   \end{equation}
   The trace of the first diagonal block can be recovered from this matrix as the sum of the first and third diagonal blocks.

   If $U_1$ and $U_2$ are chosen as the unitaries that take the nondiagonal block $(A^tA)_{12}$ to diagonal form $D_N$, things are very simplified in the above expression
   \begin{equation}
    (U_1)^t (A^tA)_{12} U_2 = D_N
   \end{equation}
   \begin{equation}\label{UnitarySubspaces}
    S^t (A^tA) S =
   \begin{pmatrix}
   (U_1)^t (A^tA)_{11} U_1 & (U_1)^t (A^tA)_{11} \bar{U}_1 & D_N & 0\\
   (\bar{U}_1)^t (A^tA)_{11} U_1 & (\bar{U}_1)^t (A^tA)_{11} \bar{U}_1 & 0 & 0 \\ 
   D_N & 0 & (U_2)^t (A^tA)_{22} U_2 & (U_2)^t (A^tA)_{22} \bar{U}_2 \\
   0 & 0 & (\bar{U}_2)^t (A^tA)_{22} U_2 & (\bar{U}_2)^t (A^tA)_{22} \bar{U}_2 
    \end{pmatrix}
   \end{equation}
 Now, a transformation should be chosen that maximises the trace of the two first blocks.   This is accomplished by a transformation that diagonalises the reduced matrix $\tilde{M}$.   The computation of such transformation does not represent a big computational cost, because of the relatively small size of this matrix.
  \begin{equation}
    \begin{pmatrix}
     (\alpha_1)^t & (\beta_1)^t\\
     (\beta_2)^t & (\alpha_2)^t
    \end{pmatrix} 
    \begin{pmatrix}
     \tilde{M}_{11} & D_N \\
     D_N & \tilde{M}_{22}
    \end{pmatrix} 
    \begin{pmatrix}
     \alpha_1 & \beta_2\\
     \beta_1 & \alpha_2
    \end{pmatrix} =
    \begin{pmatrix}
      D_1 & 0 \\
      0 & D_2
    \end{pmatrix}
  \end{equation}
  where $\tilde{M}_{11}=(U_1)^tM_{11}U_1$ and $\tilde{M}_{22}=(U_2)^tM_{22}U_2$.

  The final form of the matrix will be:
   \begin{equation}
   \tilde{U}^t S^t (A^tA) S \tilde{U} =
   \begin{pmatrix}
   D_1 & W & 0 & X \\
   W^t & \bar{\tilde{M}}_{11} & Y & 0 \\ 
   0 & Y^t & D_2 & Z \\
   X^t & 0 & Z^t & \bar{\tilde{M}}_{22} 
    \end{pmatrix}
   \end{equation}
  Where:
  \begin{multline}
   W = (\alpha_1)^t (U_1)^t(A^tA)_{11}\bar{U}_1 \hspace{2cm} X = (\beta_1)^t (U_2)^t (A^tA)_{11} \bar{U}_2\\
   Y = (\beta_2)^t  (U_1)^t(A^tA)_{11}\bar{U}_1 \hspace{2cm} Z = (\alpha_2)^t (U_2)^t (A^tA)_{11} \bar{U}_2 
  \end{multline}
  and
  \begin{equation}
   \tilde{U} = \begin{pmatrix}
    \alpha_1 & 0 & \beta_2 & 0 \\
    0 & 1 & 0 & 0\\
    \beta_1 & 0 & \alpha_2 & 0\\
    0 & 0 & 0 & 1 \end{pmatrix}
  \end{equation}
  The trace of the two first blocks will now be bigger, because the highest eigenvalues of the reduced matrix are concentrated in $D_1$.

  This procedure can be iterated, and each time the trace of the first blocks will be bigger.
  \begin{equation}
   U_{total} = S\tilde{U} = \begin{pmatrix}
    U_1 \alpha_1 & \bar{U}_1 & U_1 \beta_2 & 0\\
    U_2 \beta_1 & 0 & U_2 \alpha_2 & \bar{U}_2
   \end{pmatrix}
  \end{equation}   

 \section{An algorithm for block-SVD}
  Then, an algorithm can is proposed for block-SVD of a big matrix:
  \begin{enumerate}
   \item Make sure that there are more (or the same) rows as columns, or transpose otherwise.  This first steps are performed with the sparse csv triplet representation.
   \item Compute the rows and columns euclidean norm
   \item Order rows and columns in descending norm order
   \item Choose a cutting point for the rows and columns.  This can be made in several ways.  The one tried here is taking the point where at least 2/3 of the frobenius norm is in the first column block, and cut the row blocks as to yield a square first diagonal block.
   \item Create the blocks the appropriate size
   \item Create the Householder unitary matrix that annihilates block 12
   \item From that starting point, iterate the maximisation of the trace of block 11 of the square matrix $A^tA$ until a certain tolerance
   \item Perform SVD of the first block
   \item This gives an approximation of the eigenvalue decomposition of $A^tA$.  Multiplying the initial matrix (by blocks) times the inverse of the square root (also by blocks) the two-block relevant vertical slice of the economy unitary left ($U$ in $UDV^t$matrix are obtained.   The relevant slice of the other unitary is the first vertical 2-block slice of the unitary obtained by iteration ($V$)
  \end{enumerate}
 \section{Some results}
  A matrix of the occurrences of 3204 words in 17780 documents was used.   The procedure of cutting the blocks gave blocks with the following characteristics:\vspace{12pt}

  \begin{tabular}{|c|c|c|c|c|c|}\hline
   BLOCK & ROWS & COLUMNS & DENSITY & SQUARE NORM & NORM PERCENTAGE\\
\hline
   whole & 17780 & 3204 & 0.35\% & 5125858 & 100\% \\
   11 & 215 & 215 & 21.14\% & 3323694 & 64.84\% \\
   12 & 215 & 2989 & 3.73\% & 444596 & 8.67\% \\
   21 & 17565 & 215 & 0.51\% & 96345 & 1.88\% \\
   22 & 17565 & 2989 & 0.28\% & 1261223 & 24.60\%\\ \hline
  \end{tabular}\vspace{12pt}

  It must be noted that the matrix was cut as to leave blocks 11 and 21 with a little more than 2/3 of the square frobenius norm.

  The blocks of the matrix $A^tA$ had the following characteristics:\vspace{12pt}

  \begin{tabular}{|c|c|c|c|c|c|}\hline
   BLOCK & ROWS & COLUMNS & DENSITY & TRACE* & PERCENTAGE\\
\hline
   whole & 3204 & 3204 & 0.35\% & 5125858 & 100\% \\
   11 & 215 & 215 & 100\% & 3420039 & 66.72\% \\
   12 & 215 & 2989 & 3.73\% & 543980 & 10.61\% \\
   22 & 2989 & 2989 & 0.28\% & 1705819 & 33.28\%\\ \hline
  \end{tabular}\vspace{12pt}

  (*)The "trace" of the nondiagonal block is not actually a trace, but the sum of its singular values.   The traces of the diagonal blocks must of course sum up to the total trace, so their percentages sum up to 100\%.   The percentage for the nondiagonal block is only computed to measure how non-block-diagonal the matrix is.

  The sum of the singular values of the nondiagonal block is going to be called \textbf{nondiagonality} from now on.

  To enhance convergence, in every iteration whose number is the square of an integer, a transformation is included that tries to annihilate the block 12, but is damped by $\frac{1}{n}$, being $n$ the number of the iteration.   This is accomplished by computing matrix $\alpha$ and $\beta$ as if the whole nondiagonal block was divided by $n$.

  The values of the traces for the first five iterations are shown in the following table:\vspace{12pt}

 \begin{tabular}{|c|c|c|c|c|}\hline
  & BLOCK 11 & BLOCK 22 & BLOCK 12 & TIME (s)\\\hline
  0  &  3420039 & 1705819 & 543980 & - \\
  1  &  3864177 & 1261681 & 144997 & 4.54 \\
  2  &  3886512 & 1239346 & 24666 & 105.68\\ 
  3  &  3888110 & 1237748 & 6828 & 96.09\\
  4  &  3888370 & 1237488 & 3838 & 105.11\\
  5  &  3888479 & 1237378 & 2183 & 107.09\\
  6  &  3888541 & 1237317 & 1781 & 103.88\\
  7  &  3888584 & 1237274 & 1187 & 104.18\\
  8  &  3888617 & 1237241 & 1143 & 105.69\\
  9  &  3888644 & 1237214 &  822 & 103.51\\
  10  &  3888666 & 1237192 &  831 & 103.91\\
  ..  &  ....... &  .......  &  ........ & .....\\
  21  &  3888828 & 1237030 & 295 & 103.47\\
  22  &  3888836 & 1237021 & 269 & 102.47\\
  23  &  3888844 & 1237014 & 240 & 102.94\\
  24  &  3888850 & 1237008 & 211 & 102.77\\
  25  &  3888855 & 1237003 & 175 & 102.58\\
  26  &  3888859 & 1236999 & 158 & 101.95\\
  27  &  3888863 & 1236995 & 143 & 101.70\\
  28  &  3888866 & 1236992 & 122 & 101.78\\
  29  &  3888869 & 1236989 & 103 & 101.61\\
  30  &  3888871 & 1236987 & 91 & 101.87\\
  31  &  3888872 & 1236986 & 75 & 101.56\\ \hline
 \end{tabular}\vspace{12pt}

  The criterion used for convergence was that the ratio $\frac{Norm\ 12}{Norm\ 11}$ became 1/10000 of its initial value (about 1.16).  It can be seen that the trace of the first block is allways increased, as expected, but the sum of eigenvalues of the nondiagonal block oscillates after some iterations.

  The availability of memory for MATLAB 6.5 does not allow to perform the complete SVD decomposition of the matrix, but it is possible to compute the singular values.   The singular values contained in the first $215\times215$ block then account for 75\% of the square block of the matrix.

   The first 215 obtained singular values had differences under 1e-10 with those calculated by the usual algorithm, except for the lowest four.
   \section{Aknowledgements}
   This work was sponsored by the European Comission under the contract FP6-027026 K-Space and Foundation for the Future of Colombia COLFUTURO.    I would also like to aknowledge the valuable guidance of professor C. J. van Rijsbergen and useful advise from Mark Girolami in the developement of this work.

\end{document}